\documentclass[12pt,a4paper,twoside]{article}
\usepackage{amsfonts}
\usepackage{amsthm,amsmath,amssymb}

\def\trait #1 #2 #3 {\vrule width #1pt height #2pt depth #3pt}
\def\fin{
    \trait .3 5 0
    \trait 5 .3 0
    \kern-5pt
    \trait 5 5 -4.7
    \trait 0.3 5 0
\medskip}
\newtheorem{teor}{Theorem}[section]
\newtheorem{defin}[teor]{Definition}
\newtheorem{lemm}[teor]{Lemma}
\newtheorem{osse}[teor]{Remark}
\newtheorem{prop}[teor]{Proposition}
\newtheorem{defi}[teor]{Definition}
\newtheorem{coro}[teor]{Corollary}
\newtheorem{prob}[teor]{Problem}
\newtheorem{hypo}[teor]{Hypothesis}
\newcommand{\bele}{\begin{lemm}\begin{sl}}
\newcommand{\enle}{\end{sl}\end{lemm}}
\newcommand{\bedef}{\begin{defi}\begin{sl}}
\newcommand{\eddef}{\end{sl}\end{defi}}
\newcommand{\bete}{\begin{teor}\begin{sl}}
\newcommand{\ente}{\end{sl}\end{teor}}
\newcommand{\beos}{\begin{osse}\begin{rm}}
\newcommand{\eddos}{\end{rm}\end{osse}}
\newcommand{\bepr}{\begin{prop}\begin{sl}}
\newcommand{\empr}{\end{sl}\end{prop}}
\newcommand{\bepro}{\begin{prob}\begin{rm}}
\newcommand{\empro}{\end{rm}\end{prob}}
\newcommand{\bede}{\begin{defin}\begin{sl}}
\newcommand{\edde}{\end{sl}\end{defin}}
\newcommand{\beco}{\begin{coro}\begin{sl}}
\newcommand{\enco}{\end{sl}\end{coro}}
\newcommand{\behy}{\begin{hypo}\begin{sl}}
\newcommand{\enhy}{\end{sl}\end{hypo}}
\newcommand{\prova}{\noindent{\bf Proof.\hspace{4mm}}}
\newcommand{\thspace}{\hspace{3mm}}

\newcommand{\II}{\mathbb{I}}
\newcommand{\RR}{\mathbb{R}}
\newcommand{\NN}{\mathbb{N}}

\newcommand{\beeq}[1]{\begin{equation}\label{#1}}
\newcommand{\eddeq}{\end{equation}}
\newcommand{\beeqa}[1]{\begin{eqnarray}\label{#1}}
\newcommand{\eddeqa}{\end{eqnarray}}
\newcommand{\beal}[1]{\begin{align}\label{#1}}
\newcommand{\eddal}{\end{align}}
\newcommand{\bespl}[1]{\begin{split}\label{#1}}
\newcommand{\edspl}{\end{split}}
\newcommand{\bega}[1]{\begin{gather}\label{#1}}
\newcommand{\edga}{\end{gather}}
\newcommand{\beeqax}{\begin{eqnarray*}}
\newcommand{\eddeqax}{\end{eqnarray*}}
\def\qed{\ifmmode   \else \leavevmode\unskip\penalty9999 \hbox{}\nobreak\hfill
  \fi
  \quad\hbox{\hskip.5em\vrule width.4em height.6em depth.05em\hskip.1em}}
\def\endproofsym{\qed}

\def\endnobox{\def\endproofsym{}\end{proof}\def\endproofsym{\qed}}

\newcommand{\no}{\nonumber}
\newcommand{\beeqao}{\begin{eqnarray}\no}
\newcommand{\bealo}{\begin{align}\no}
\newcommand{\besplo}{\begin{split}\no}
\newcommand{\begao}{\begin{gather}\no}
\def\trait #1 #2 #3 {\vrule width #1pt height #2pt depth #3pt}
\def\fin{\hfill
    \trait .3 5 0
    \trait 5 .3 0
    \kern-5pt
    \trait 5 5 -4.7
    \trait 0.3 5 0
\medskip}

\newcommand{\duav}[1]{\langle{#1}\rangle}

\newcommand{\duaw}[1]{_{W'}\langle{#1}\rangle_{W}}

\newcommand{\dt}{\partial_t}

\newcommand{\itt}{\int_0^t}

\newcommand{\io}{\int_\Omega}

\newcommand{\e}{\varepsilon}
\newcommand{\ee}{^{\varepsilon}}

\newcommand{\mezzo}{{\frac{1}{2}}}

 \DeclareMathOperator{\dive}{div}

\newcommand{\HUH}{W^{1,2}(0,T;H)}
\newcommand{\HUV}{W^{1,2}(0,T;V)}

\newcommand{\CZH}{C^0([0,T];H)}
\newcommand{\CZV}{C^0([0,T];V)}

\newcommand{\LDH}{L^2(0,T;H)}
\newcommand{\LDV}{L^2(0,T;V)}
\newcommand{\LDVp}{L^2(0,T;V')}

\newcommand{\LIV}{L^\infty(0,T;V)}

\newcommand{\LIW}{L^{\infty}(0,T;W)}

\newcommand{\HUtV}{W^{1,2}(0,t;V)}
\newcommand{\HUtW}{W^{1,2}(0,t;W)}

\newcommand{\CZtH}{C^0([0,t];H)}
\newcommand{\CZtV}{C^0([0,t];V)}

\newcommand{\LDtH}{L^2(0,t;H)}
\newcommand{\LDtV}{L^2(0,t;V)}
\newcommand{\LDtVp}{L^2(0,t;V')}
\newcommand{\LItH}{L^\infty(0,t;H)}

\let\TeXchi\chi
\def\chi{{\setbox0 \hbox{\mathsurround0pt
$\TeXchi$}\hbox{\raise\dp0 \copy0 }}}

\newcommand{\betapb}{\bar{\beta_1}}
\newcommand{\betasb}{\bar{\beta_2}}
\newcommand{\betatb}{\bar{\beta_3}}

\newcommand{\xib}{\boldsymbol{\xi}}

\newcommand{\chib}{{\boldsymbol{\chi}}}
\newcommand{\Ub}{{\bf U}}

\newcommand{\betab}{{\boldsymbol{\beta}}}
\newcommand{\ub}{{\bf u}}
\newcommand{\vb}{{\bf v}}

\newcommand{\teta}{\vartheta}

\def\fine{\hfill\kern4pt \vrule height4pt depth0pt width4pt }

\def\dive{\mbox{\rm div\,}}

\numberwithin{equation}{section}
\numberwithin{equation}{section}
\begin{document}

\title{A model for shape memory alloys\\
with the possibility of voids}

\date{}


\author{\renewcommand{\thefootnote}{\!$\fnsymbol{footnote}$}
Michel Fr\'emond\\
{\sl  Dipartimento di Ingegneria Civile, Universit\`a di Roma ``Tor Vergata''}\\
 {\sl Via del Politecnico, 1, I-00133 Roma, Italy}\\
{\rm E-mail:~~\tt Michel.Fremond@uniroma2.it}
\vspace{0.4truecm}
\and
Elisabetta Rocca\\
{\sl Dipartimento di Matematica, Universit\`a di Milano,}\\
{\sl Via Saldini 50, 20133 Milano, Italy}\\
{\rm E-mail:~~\tt elisabetta.rocca@unimi.it} } \maketitle

Dedicated to Professor Roger Temam on the occasion of his 70th
anniversary
\bigskip

\begin{abstract}
The paper is devoted to the study of a mathematical model for the
thermomechanical evolution of metallic shape memory alloys. The
main novelty of our approach consists in the fact that we include
the possibility for these materials to exhibit voids during the
phase change process. Indeed, in the engineering paper
\cite{allvoids} has been recently proved that voids may appear
when the mixture is produced by the aggregations of powder. Hence,
the composition of the mixture varies (under either thermal or
mechanical actions) in this way: the martensites and the austenite
transform into one another whereas the voids volume fraction
evolves. The first goal
of this contribution is hence to state a PDE system capturing all these modelling
aspects in order then to establish the well-posedness of the
associated initial-boundary value problem.

\end{abstract}
\vspace{.4cm}

\noindent {\bf Key words:}\thspace shape memory alloys, mixtures
with voids,
 nonlinear PDEs system, initial-boundary
value problem, existence of solutions, continuous dependence
result

\vspace{4mm}

\noindent
{\bf AMS (MOS) subject clas\-si\-fi\-ca\-tion:}\thspace 80A22, 74D10, 35A05, 35Q30, 35G25.

\pagestyle{empty}



\pagestyle{myheadings}
\newcommand\testopari{\sc Fr\'emond--Rocca}
\newcommand\testodispari{\sc SMA with voids}
\markboth{\testodispari}{\testopari}



\section{Introduction}
\label{sec:intro}

Shape memory alloys are mixtures of many martensites variants and
of austenite. They exhibit an unusual behavior: even if they are
permanently deformed, they can totally recover their initial shape
just by thermal or mechanical means. There may be voids in the
mixture, which may appear when the mixture is produced by the
aggregations of powders, as it has been recently proved in the
engineering paper \cite{allvoids}. The composition of the mixture
varies: the martensites and the austenite transform into one
another whereas the voids volume fraction evolves. These phase
changes can be produced either by thermal actions or by mechanical
actions. The striking properties of shape memory alloys result
from interactions between mechanical and thermal actions (cf.,
e.g., \cite{BerrietLexcellent1992, Guenin1986}).

We assume that the phases can coexist at each point and we suppose that,
besides austenite, only two martensitic variants are present. However, this choice
provides a sufficiently good description of the phenomenon, as we want
describe a macroscopic predictive theory which can be used for
engineering purposes. The phase volume fractions, which are state
quantities, are subjected to constraints. In particular,  their
sum must be lower than $1$, but not necessarily equal to 1, due to
the presence of voids (cf.~also \cite{frm3as} and \cite{frquart}
where this property is introduced in order to treat solid-liquid
phase transitions with the possibility of voids and \cite{famf}
for the corresponding numerical results). It is shown that most of
the properties of shape memory alloys result from careful
treatment of those internal constraints (cf., e.g.,
\cite{Fremond1987b}--\cite{FremondMyasaki1996}).
All these quoted references are related to a three dimensional model taking the temperature,
the macroscopic deformation and the volumetric proportion of austenite and martensite
as state variables.
Moreover, let us note that there are not too many mathematical models
describing phase transitions in which the interactions between
different types of substances and the possibility of having voids
is taken into account: we can quote only the two contributions
\cite{frm3as,frquart}.

It is beyond our purposes to give a complete description of the
existing literature on models for SMA. However, restricting
ourselves to the macroscopic description of these phenomena, we
can refer to the main contributions
\cite{Fremond1987b}--\cite{FremondMyasaki1996}, \cite{ap1, ap2,
falk, fk, ps06} and \cite{ams, bonetti1, bonetti2,
ChrysochoosPham1993, ChrysochoosLobel1994, lobel1994, Pagano1997,
ypz} (and references therein) describing full thermomechanical
models and studying the resulting PDEs from mathematical viewpoint
respectively. We shall instead focus here on a generalization of
the Fr\'emond model for SMA introduced in
\cite{Fremond1987b}--\cite{FremondMyasaki1996}.

Let us then explain in detail which is the main aim of this contribution, compare it with the
results already present in the literature, and show the main
mathematical difficulties encountered. As already mentioned, in this paper we deal
with a generalization of the model introduced in \cite{Fremond1987b}--\cite{FremondMyasaki1996}
and later on studied in many contributions starting from the pioneering paper
\cite{cfv}, where  an existence and uniqueness result has been proved for the solution of a
simplified problem, where all the nonlinearities in the balance of energy
are neglected and the momentum balance equation is considered in the quasi-stationary
form and fourth order terms (related to the second gradient theory) are taken into account.
In the case when the fourth-order term is omitted an existence result dealing with the
linearized energy balance equation has
been proved in \cite{colli}, while \cite{collibis} one can find the proof of the existence of
solutions to the linearized problem by including an inertial term in the momentum balance.
We can report also of some results when some or all the nonlinearities are kept in the energy balance.
The full one-dimensional model is shown to admit a unique solution both in the quasi-stationary
case in \cite{ColliSprekels1995} and in the case of a hyperbolic momentum equation in
\cite{ColliSprekels1996, smetov}. Existence
results have been proved also for the three-dimensional model (cf. \cite{colliter, ColliSprekels1992, hnz}).
Finally, let us mention the  uniqueness result for the full quasi-static three-dimensional model proved in \cite{Chemetov1997}
and an updated and detailed presentation of the Fr\'emond model and related system of
equations and conditions, applying to the multidimensional case as well, which is provided in \cite{bonetti1, bonetti2},
\cite[Chapter 13]{Fremond2001}, and \cite{FremondMyasaki1996}. Let us also point out \cite{bonetti1, bonetti2} for recent existence and
uniqueness results in the three-dimensional situation, where the various nonlinear terms
arising in the derivation of the model are accounted. The
large time behavior of solutions is investigated in \cite{cls} in connection with the convergence
to steady-state solutions and in \cite{cshirak, cfrs} where the authors characterize the large time behavior according to the
theory of dissipative dynamical systems.

However, all these contributions were dealing with the case in
which no voids can occur between phases. To model this possibility
and to solve rigorously the results PDE system is just our aim
here. First, in the next Section~\ref{sec:model}, we derive a
model taking the possibility of having voids into account,
introducing rigorously the pressure which has a paramount
importance on the mechanical behaviour. In order to do that, we
follow the ideas of \cite{frm3as} in which this was done in case
of a two-phase transition phenomenon. Then, in
Section~\ref{mainres}, we give a rigorous formulation of initial
boundary value problem associated with the resulting PDEs and we
state our main results: existence, uniqueness, and continuous
dependence of solutions from the data. The proofs are carried over
in Section~\ref{exiproof} and Section~\ref{prouni}. The main
mathematical difficulties are due to the nonlinear and singular
coupling between the equations. In particular, in order to
describe the evolution of the absolute temperature variable, we
shall use the entropy balance equation (cf.~\cite{bcf}--\cite{bfr}
for a complete derivation and motivation of this equation). This
equation turns out to be singular in $\teta$ but the main
advantage of using it is that once one has proved that a solution
component $\teta$  does exists then it turns automatically out to
be positive and the proof of positivity of the absolute
temperature is historically one of the main difficulties of these
types of problems. The idea here is to approximate the
nonlinearities with regular functions, to solve the regularized
system by means of a Banach fixed point argument and then to use
compactness and lower semicontinuity arguments in order to pass to
the limit and obtain a solution of the original problem.


\section{The derivation of the model}
\label{sec:model}

In this section we explicitly derive a macroscopic model
describing the evolution od SMA with the possibility of voids. The
model is obtained by properly choosing the state quantities, the
balance laws and the constitutive relations in agreement with the
principle of thermodynamics and with experimental evidence.

\subsection{The State Quantities}
\label{subsec:state}

We deal only with macroscopic phenomena and macroscopic quantities. To
describe the deformations of the alloy, the macroscopic small deformation $%
\varepsilon (\ub)$, ($\ub$ being the small displacement) and the
temperature $\teta$ are chosen as state quantities.

The properties of shape memory alloys result from martensite--austenite
phase changes produced either by thermal actions (as usual) or by
mechanical actions. On the macroscopic level, some quantities are needed to take
those phase changes into account. For this purpose, the volume fractions $%
\beta _{i}$ of the martensites and austenite are chosen as state quantities.
For simplicity, we assume that only two martensites exist together with
austenite. The volume fractions of the martensites are $\beta _{1}$ and $%
\beta _{2}$. The volume fraction of austenite is $\beta _{3}$. These volume
fractions are not independent: they satisfy the following {\sl internal} constraints
\begin{equation}
0\leq \beta _{i}\leq 1, \quad i=1,2,3 \label{1}
\end{equation}%
due to the definition of volume fractions. Since we assume that voids can appear
in the martensite--austenite mixture, then the $\beta $'s must satisfy an other
internal constraint
\begin{equation}
\beta _{1}+\beta _{2}+\beta _{3}\leq 1 \label{2}
\end{equation}%
the quantity $v=1-\left( \beta _{1}+\beta _{2}+\beta _{3}\right) $
being the
voids
 volume fraction. This is the case when the alloy is produced
by aggregating powders as shown in \cite{allvoids}. In case no
voids are considered, this sum should be equal to 1 and this
considerably simplifies the analysis (cf., e.g., \cite{cfv}).

We denote by $\betab$\textbf{\ }the vector of components $\beta _{i}$ ($i=1,2,3$) and the set
of the state quantities is
\begin{equation*}
E=\left\{ \varepsilon(\ub) ,\betab,\nabla\betab,\teta\right\}
\end{equation*}%
while the quantities which describe the evolution and the thermal heterogeneity are%
\begin{equation*}
\delta E=\{\varepsilon(\ub)_t,\betab_t,\nabla\betab_t,\nabla \teta\}\,.
\end{equation*}
The gradient of $\nabla\betab$ accounts for local interactions of
the volume fractions at their neighborhood points.

\subsection{The mass balance}
\label{subsec:mb}

Assuming the same constant density $\rho $ (the reader can refer to \cite{frquart} for a model in which
different densities of the substances are taken into account in a general two-phase change phenomenon)
and the same velocity $\Ub=\ub_t$ for each phase, the mass balance reads
\begin{equation*}
\rho(\beta _{1}+\beta _{2}+\beta _{3})_t+\rho (\beta _{1}+\beta
_{2}+\beta _{3}){\dive}\Ub=0\,.
\end{equation*}%
Within the small perturbation assumption, this equation gives%
\begin{equation*}
\rho (\beta _{1}+\beta _{2}+\beta _{3})_t+\rho
(\beta _{1}^{0}+\beta _{2}^{0}+\beta _{3}^{0}){\dive}\Ub=0
\end{equation*}%
where the $\beta _{i}^{0}$'s are the initial values of the $\beta _{i}$. For
the sake of simplicity, we assume
\begin{equation}
\beta _{1}^{0}+\beta _{2}^{0}+\beta _{3}^{0}=1\label{1000}
\end{equation}%
and have $\dt(\beta _{1}+\beta _{2}+\beta _{3})+{\dive}\Ub=0$, hence
\begin{gather}
\dt(\beta _{1}+\beta _{2}+\beta _{3})+\dive\ub_t=0\,.  \label{1002}
\end{gather}%
Mass balance is a relationship between the quantities of $\delta
E$, indeed, its effects will be included in the cinematic
 relations (cf.~\eqref{fi} in the following subsections).

\subsection{The equations of motion}
\label{subsec:motion}

They result from the principle of virtual power involving the
power of the internal forces, (cf, e.g., \cite{Fremond2001})%
\begin{equation*}
-\int_{\Omega }\left\{ \sigma :D({\bf V})+{\bf B}\cdot {\bf
\delta}+H:\nabla{\bf \delta}\right\} d\Omega
\end{equation*}%
where ${\bf V}$ and ${\bf \delta}$ are virtual velocities, the actual
velocities being $\Ub$ and $\betab_t$. The internal forces are the
stress $\sigma $, the phase change work
vector ${\bf B}$, and the phase change work flux tensor $H$. The
equations of motion are
\begin{gather}
\rho \Ub_t=\dive\sigma+{\bf f},\quad
0={\dive}H-{\bf B}+{\bf A}\quad\text{in}\;\Omega  \label{3} \\
\sigma {\bf n}={\bf g},\;H{\bf n}={\bf a}\quad \text{on}\;\partial \Omega
\label{3a}
\end{gather}%
where $\rho $ is the density, $\Ub_t$ the acceleration of the alloy
which occupies the domain $\Omega $, with boundary $\partial \Omega $ and
outward normal vector ${\bf n}$. The alloy is loaded by body forces ${\bf f}$
and by surface tractions ${\bf g}$, and submitted to body sources of damages
${\bf A}$ and surfaces sources of damage ${\bf a}$ (in the following we will suppose,
for simplicity, ${\bf A}={\bf a}={\bf 0}$).

\subsection{The free energy}
\label{subsec:free}

As explained above, a shape memory alloy is considered as a mixture of the
martensite and austenite phases with volume fractions $\beta _{i}$. The
volume free energy of the mixture we choose is
\begin{equation}
\Psi =\Psi(E)=\sum\limits_{i=1}^{3}\beta _{i}\Psi _{i}(E)+h(E) \label{10}
\end{equation}%
where the $\Psi _{i}$'s are the volume free energies of the ${i}$ phases
and $h$ is a free energy describing interactions between the different
phases. We have assumed that internal constraints are physical properties, hence,
we decide to choose properly  the two
functions describing the material, i.e., the free energy $\Psi $
and the pseudopotential of dissipation $\Phi $, in order to take these constraints
into account. Since, the pseudopotential
describes the kinematic properties (i.e., properties which depend on the
velocities) and the free energy describes the state properties, obviously the
internal constraints (\ref{1}) and (\ref{2}) are to be taken into account with the choice of
the free energy $\Psi $.

For this purpose, we assume the $\Psi _{i}$'s are defined over the whole
linear space spanned by $\beta _{i}$ and the free energy is defined by
\begin{equation*}
\Psi (E)=\beta _{1}\Psi _{1}(E)+\beta _{2}\Psi _{2}(E)+\beta _{3}\Psi
_{3}(E)+h(E)\,.
\end{equation*}%
We choose the very simple interaction free energy%
\begin{equation*}
h(E)=I_{C}(\betab)+\frac{k}{2}\left|\nabla\betab\right|^{2}
\end{equation*}%
where $I_{C}$ is the indicator function of the convex set
\begin{equation}\label{defC}
C=\{(\gamma _{1},\gamma _{2},\gamma _{3})\in
\mathbb{R}^{3};0\leq \gamma _{i}\leq 1;\gamma _{1}+\gamma _{2}+\gamma _{3}\leq 1\}\,.
\end{equation}
Moreover, and by $(k/2)|\nabla \betab|^{2}$ we mean the product of
two tensors $ \nabla\betab$ multiplied by the \textsl{interfacial
energy coefficient} $(k/2)>0$. The terms
$I_{C}(\betab)+(k/2)\vert\nabla\betab\vert^2$ may be seen as a
\textsl{ mixture or interaction free-energy}.

The only effect of $I_{C}(\betab)$ is to guarantee that the
proportions $\beta _{1}$,
$\beta _{2}$ and $\beta _{3}$
take admissible physical values, i.e.~they satisfy constraints
(\ref{1}) and (\ref{2}) (cf.~also \eqref{defC}). The interaction
free energy term $I_{C}(\betab)$ is equal to zero when the mixture
is physically possible ($\betab\in C$) and to $+\infty $ when the
mixture is physically impossible ($\betab\notin C$).

Let us note even if the free energy of the voids phase is $0$, the
voids phase has physical properties due to the interaction free
energy term $(\nu/2)\vert\nabla\betab\vert^2$ which depends on the
gradient of $\betab$. It is known that this gradient is related to
the interfaces properties: $ \nabla {\beta _{1}}$, $ \nabla {\beta
_{2}}$ describes properties of the voids-martensites interfaces
and $ \nabla {\beta _{3}}$ describes properties of the
voids-austenite interface. In this setting, the voids have a role
in the phase change and make it different from a phase change
without voids. The model is simple and schematic but it may be
upgraded by introducing sophisticated interaction free energy
depending on $\betab$ and on $\nabla\betab$.

For the volume free energies, we choose
\begin{gather*}
\Psi _{1}(E)=\frac{1}{2}\varepsilon(\ub)
:K_{1}:\varepsilon(\ub) +\sigma _{1}(\teta):\varepsilon(\ub) -C_{1}\teta\log \teta,
\\
\Psi _{2}(E)=\frac{1}{2}\varepsilon(\ub)
:K_{2}:\varepsilon(\ub) +\sigma _{2}(\teta):\varepsilon(\ub) -C_{2}\teta\log \teta,
\\
\Psi _{3}(E)=\frac{1}{2}\varepsilon(\ub)
:K_{3}:\varepsilon(\ub) -\frac{l_a}{\teta_{0}}(\teta-\teta_{0})-C_{3}\teta\log \teta,
\end{gather*}%
where $K_{i}$ are the volume elastic tensors and $C_{i}$ the
volume heat capacities of the phases. Stresses $\sigma _{i}(\teta)$ depend on
temperature $\teta$ and the quantity $l_{a}$ is the latent heat
martensite-austenite volume phase change at temperature $\teta_{0}$ (see Remark~\ref{teta} below).

\beos\label{teta}
To make the model more realistic, we can introduce two temperatures to
characterize the transformation: $\teta_{0}$, the temperature at the beginning
of the transformation and $\teta_{f}$ the temperature at the end. The
interaction free energy is completed by $h(\betab)=(l_{a}/\teta_{0})
(\teta_{0}-\teta_{f})(\beta _{3})^{2}$ (cf.~\cite{BalandraudErnest1998}, \cite{ps06}, \cite{Pham1994}).
\eddos

Because we want to describe the main basic properties of the shape memory
alloys with voids, we assume that the elastic matrices $K_{i}$ and the heat
capacities $C_{i}$ are the same for all of the phases:
\begin{equation*}
C_{i}=\bar C,\;K_{i}=K\quad i=1,2,3\,.
\end{equation*}%
Always for the sake of simplicity, we assume that
\begin{equation*}
\sigma _{1}(\teta)=-\sigma _{2}(\teta)=-\tau (\teta)\II
\end{equation*}%
where $\II$ stands for the
identity
 matrix. Concerning the stress $\tau (\teta)$, it is known
that at high temperature the alloy has a classical elastic
behaviour. Thus $\tau (\teta)=0$ at high temperature, and we
choose the schematic simple expression
\begin{equation*}
\tau (\teta)=(\teta-\teta_{c})\overline{\tau },\;\text{for\ }\teta\leq \teta_{c},\tau (\teta)=0,\;%
\text{for\ }\teta\geq \teta_{c}
\end{equation*}%
with $\overline{\tau }\leq 0$ and assume the temperature $\teta_{c}$ is
greater than $\teta_{0}$. With those assumptions, it results
\begin{gather*}
\Psi (E)=\frac{(\beta _{1}+\beta _{2}+\beta _{3})}{2}\left\{ \varepsilon(\ub)
:K:\varepsilon(\ub) \right\} \\
-(\beta _{1}-\beta _{2})\tau (\teta)\II:\varepsilon(\ub) -\beta _{3}\frac{l_{a}}{\teta_{0}}%
(\teta-\teta_{0})-C\teta\log \teta+\frac{k}{2}\left|\nabla\betab\right|^{2}+I_{C}(\betab)\,.
\end{gather*}

\subsection{The pseudo-potential of dissipation}
\label{subsec:pseudo}

The dissipative forces are defined via a pseudo-potential of dissipation $%
\Phi $ introduced by J.J. Moreau (it is a convex, positive
function with value zero at the origin, \cite{et}, \cite{mt},
\cite{moreau 1966}). As already
remarked, the mass balance (\ref{1002}) is a relationship between velocities of $%
\delta E$. Thus we take it into account in order to define the
pseudo-potential and introduce the indicator function $I_{0}$ of
the origin of
$\mathbb{R}$
 as follows
\begin{equation*}
I_{0}(\dt(\beta _{1}+\beta _{2}+\beta _{3})+\dive \ub_t)\,.
\end{equation*}
>From experiments, it is known that the behaviour of shape memory alloys
depends on time, i.e., the behaviour is dissipative. We define a
pseudopotential of dissipation
\begin{equation}\label{fi}
\Phi (\teta,\betab_t,\nabla \teta)=%
\frac{c}{2}\left|\betab_t\right|
^{2}+\frac{\upsilon}{2}|\nabla\betab_t|^{2}+\frac{\lambda }{2\teta}\left|\nabla \teta\right|^{2}+I_{0}(\dt(\beta
_{1}+\beta _{2}+\beta _{3})+\dive \ub_t)
\end{equation}
where $\lambda \geq 0$ represents the thermal conductivity and
 $c\geq 0$, $\upsilon\geq 0$ stand for phase
change viscosities.

\subsection{The constitutive laws}
\label{subsec:consti}

The internal forces are split between non-dissipative forces $\sigma ^{nd}$,
${\bf B}^{nd}$ and $H^{nd}$ depending on $(E,x,t)$ and
dissipative forces by $\left\{ \sigma ^{d},{\bf B}^{d},H^{d},-%
{\bf Q}^{d}\right\} $ depending on $\delta E=\{\varepsilon(\ub)_t,
\betab_t,\nabla\betab_t,\nabla \teta\}$ and $(E,x,t)$

\begin{equation*}
\sigma =\sigma ^{nd}+\sigma ^{d},\ \ {\bf B}%
={\bf B}^{nd}+{\bf B}^{d},\ H=H^{nd}+H^{d}
\end{equation*}%
with the entropy flux vector ${\bf Q}$ being%
\begin{equation*}
{\bf Q}={\bf Q}^{d}\,.
\end{equation*}%
The nondissipative forces are defined with the free energy
\begin{gather}
\sigma ^{nd}(E)=\frac{\partial \Psi }{\partial \varepsilon(\ub) }(E)=(\beta
_{1}+\beta _{2}+\beta _{3})K:\varepsilon(\ub) -(\beta _{1}-\beta _{2})\tau (\teta)\II
\label{13} \\
{\bf B}^{nd}(E)=\frac{\partial \Psi }{\partial \betab}%
(E)=\frac12\left(
\begin{array}{c}
\varepsilon(\ub) :K:\varepsilon(\ub)  -2\tau (\teta):\varepsilon(\ub) \\
\varepsilon(\ub) :K:\varepsilon(\ub) +2\tau (\teta):\varepsilon(\ub) \\
\varepsilon(\ub) :K:\varepsilon(\ub)
-2\frac{l_{a}}{\teta_{0}}(\teta-\teta_{0})
\end{array}%
\right)
+{\bf B}^{ndr}(E) \label{13b} \\
{\bf B}^{ndr}(E)\in \partial I_{C}(\betab)
\label{13c} \\
H^{nd}=k\nabla\betab \label{15}
\end{gather}%
and the dissipative forces are defined with the pseudo-potential of
dissipation%
\begin{equation}
\left\{ \sigma ^{d},{\bf B}^{d},H^{d},-{\bf Q}^{d}\right\}
=\partial \Phi (E,\delta E)  \label{21}
\end{equation}%
where the subdifferential of $\Phi $ is with respect to $\delta E$.
Relationship (\ref{21}) gives%
\begin{gather}
\sigma ^{d}=-p\II \label{1003} \\
{\bf B}^{d}=-p\left(
\begin{array}{c}
1 \\
1 \\
1%
\end{array}%
\right) +c\betab_t  \label{1005} \\
H^{d}=\upsilon\nabla\betab_t  \label{1006} \\
-{\bf Q}^{d}=\frac{\lambda }{\teta}\nabla \teta  \label{1007}
\end{gather}%
where $p$ is the pressure in the mixture and it results
\begin{equation}
-p\in \partial I_{0}(\dt(\beta _{1}+\beta _{2}+\beta
_{3})+\dive\ub_t)\,.  \label{1001}
\end{equation}
The state laws (\ref{13})--(\ref{15}), besides implying that the
internal constraints are satisfied, give also the value of the
reactions, during the evolution, to these internal constraints.

Relationships (\ref{13})--(\ref{15}) and
(\ref{1003})--(\ref{1001}) give the constitutive laws
\begin{gather}
\sigma =(\beta _{1}+\beta _{2}+\beta _{3})K:\varepsilon(\ub) -\left((\beta _{1}-\beta
_{2})\tau (\teta)+p\right)\II,  \label{22} \\
{\bf B}(E,\delta E)=\frac12\left\{
\begin{array}{c}
\varepsilon(\ub) :K:\varepsilon(\ub) -2\tau (\teta):\varepsilon(\ub) -p\\
\varepsilon(\ub) :K:\varepsilon(\ub) +2\tau (\teta):\varepsilon(\ub) -p\\
\displaystyle\varepsilon(\ub) :K:\varepsilon(\ub)
-2\frac{l_{a}}{\teta_{0}}(\teta-\teta_{0})-p
\end{array}
\right\} +{\bf B}^{ndr}(E)+c\betab_t,  \label{23} \\
{\bf B}^{ndr}(E)\in \partial I_{C}(\betab)
\label{23ab} \\
-p\in \partial I_{0}(\dt(\beta _{1}+\beta _{2}+\beta _{3})+\dive\ub_t) \label{23ac} \\
H=k\nabla\betab+\upsilon\nabla\betab_t  \label{23a} \\
-{\bf Q}(E,\delta E)=-{\bf Q}^{d}(E,\delta E)=-\frac{\lambda }{\teta}\nabla%
\teta\,.  \label{25}
\end{gather}
It can be proved that our choice is such that the internal
constraints and the second law of thermodynamics are satisfied
(cf., e.g., \cite{FremondMyasaki1996,Fremond2001} and the next
Subsection~\ref{subsec:entrbal}).

\subsection{The entropy balance}
\label{subsec:entrbal}

By denoting
\begin{equation}
s=-\frac{\partial \Psi }{\partial \teta}=\bar C(1+\log \teta)+\beta _{3}\frac{l_{a}}{\teta_{0}}\label{1008}
\end{equation}
the entropy balance is
\begin{gather}
\frac{ds}{dt}+\dive{\bf Q}=R+\frac{1}{\teta}\left\{ \sigma ^{d}:\varepsilon (\ub)_t
+{\bf B}^{d}\frac{\partial \betab}{\partial t}+H^{d}:
\nabla\betab_t-{\bf Q}\cdot \nabla
\teta\right\}  \notag \\
=R+\frac{1}{\teta}\left\{ c|\betab_t|^{2}+\upsilon |\nabla\betab_t|^2
+\frac{\lambda }{\teta}|\nabla \teta|^{2}\right\} ,\;\text{in\ }\Omega
\label{7a} \\
-{\bf Q}\cdot {\bf n}=\Pi ,\;\text{in\ }\Omega  \label{8}
\end{gather}%
because
\begin{equation*}
p\left((\beta _{1}+\beta _{2}+\beta _{3})_t +{\dive}\Ub\right)=0
\end{equation*}%
due to (\ref{1001}), $\teta{\bf Q}$ is the heat flux vector,
$R\teta$ is the exterior volume rate of heat that is supplied to
the alloy, $\teta\pi $ is the rate of heat that is supplied by
contact action, $\varepsilon(\ub_t)$ is the strain rate. The
constitutive laws, within the small perturbation assumption and
(\ref{1000}), become
\begin{gather}
\sigma =K:\varepsilon(\ub) -\left((\beta _{1}-\beta _{2})\tau (\teta)+p\right)\II
\label{1009} \\
-p\in \partial I_{0}((\beta _{1}+\beta _{2}+\beta _{3})_t+{\dive}\Ub)\\
{\bf B}=\left\{
\begin{array}{c}
-\tau (\teta):\varepsilon(\ub) -p \\
\tau (\teta):\varepsilon(\ub) -p \\
\displaystyle-\frac{l_{a}}{\teta_{0}}(\teta-\teta_{0})-p%
\end{array}%
\right\} +{\bf B}^{ndr}+c\betab_t \\
{\bf B}^{ndr}\in \partial I_{C}(\betab) \\
H=k\nabla\betab+\upsilon\nabla\betab_t \\
{\bf Q}(E,\delta E)={\bf Q}^{d}(E,\delta E)=-\frac{\lambda }{\teta}\nabla \teta\,.
\label{1010}
\end{gather}

\subsection{The set of partial differential equations}

We assume also quasi-static evolution and, using again the small perturbation
assumption, we get the following set of partial
differential equations coupling the equations of motion (\ref{3}), the entropy
balance (\ref{7a}) and constitutive laws (\ref{1009})--(\ref{1010})%
\begin{gather}\label{pde1}
{\dive}\left((K:\varepsilon(\ub) -\left((\beta _{1}-\beta _{2})\tau (\teta)+
p\right)\II)\right)+{\bf f}=0\\
\label{pde2}
-p\in \partial I_{0}(\dt(\beta _{1}+\beta _{2}+\beta _{3})+\dive\ub_t) \\
\label{pde3}
c\betab_t-\upsilon\Delta\betab_t-k\Delta \betab+\left(
\begin{array}{c}
-\tau (\teta):\varepsilon(\ub) -p \\
\tau (\teta):\varepsilon(\ub) -p \\
-\frac{l_{a}}{\teta_{0}}(\teta-\teta_{0})-p%
\end{array}
\right)+{\bf B}^{ndr}=0 \\
\label{pde4}
{\bf B}^{ndr}\in \partial I_{C}(\betab)\\
\label{pde5}
\bar C\frac{\partial \log \teta}{\partial t}+\frac{l_{a}}{\teta_{0}}\dt\beta
_{3}-\lambda \Delta \log \teta=R.
\end{gather}
This set is completed by suitable initial conditions and the following boundary conditions:
\begin{align}\label{bousigma}
& \sigma \mathbf{n}=\mathbf{g}\quad \hbox{on }\Sigma _{1}:=\Gamma_{1}\times \lbrack 0,T]\\
\label{bouu}
& \mathbf{u}=\mathbf{u}_{t}=0\quad \hbox{on }\Sigma _{0}:=\Gamma_{0}\times \lbrack 0,T]\\
\label{boubeta}
& \partial _{\mathbf{n}}{\dot{{\boldsymbol{\beta}}}}+\partial _{\mathbf{n}}{
\boldsymbol{\beta}}=0\quad \hbox{on }\Sigma:=\partial
\Omega \times \lbrack 0,T]\\
\label{boulog}
&\partial _{\mathbf{n}}(\ln \vartheta )=\Pi \quad \hbox{on }\Sigma
\end{align}%
where $\partial_\mathbf{n}$ is the normal outward derivative to the surface $\partial \Omega $, $\mathbf{g}$ is the exterior
contact force applied to $\Gamma_{1}$, where ($\Gamma_0$,\,$\Gamma_1$) is a partition of $\partial\Omega$ and $
\Gamma_0$, $\Gamma_{1}$ have positive measures.

\subsection{Remarks on the model}

The evolution of a structure made of shape memory alloys, i.e., the
computation of $E(x,t)=(\varepsilon(\ub) (x,t),\beta _{1}(x%
,t),\beta _{2}(x,t),\beta _{3}(x,t),\teta(x,t))$ depending on
the point $x$ of the domain $\Omega$ occupied by the structure and on time $t$,
can be performed by solving numerically the set of partial differential
equations resulting from the equations of motion (\ref{3}), (\ref{3a}) the
energy balance (\ref{7a}), (\ref{8}) and the constitutive laws (\ref{22})--(%
\ref{25}), completed by convenient initial and boundary conditions
(cf., e.g., \cite{ColliSprekels1995}, \cite{Chemetov1997}, \cite{Worshing1995},
\cite{NiesgodskaSprekels1991}). The model we have described here is able to account
for the different features of the shape memory alloys: in particular, their macroscopic, mechanical
and thermal properties. We have used schematic free energies and schematic
pseudopotentials of dissipation.

There are still many possibilities to upgrade the basic choices we have made to
take into account the practical properties of shape memory alloys. Let us, for instance,
mention that the pseudopotential of dissipation can be
modified in order to describe more precisely the hysteretical properties of the materials. There is no
difficulty in having more than two martensites, for instance, to take care
of $24$ possible martensites! In the same way, it is possible to
take into account of the different forms of a single martensite variant, as explained in
\cite{ps06}.

Note that the physical quantities for characterizing an educated shape
memory alloys are $K$, $\bar C$, $l_{a}$, $\teta_{0}$, $\teta_{c}$, $\tau $, the two
martensite volume fractions (for the free energy) and $c$, $k$, $\lambda $ (for the pseudopotential
of dissipation). They  are  indeed not so many in order to have a complete multidimensional
model which can be used for engineering purposes.

Other models and results may be found in \cite{ap1, ap2, BerveillerPatoor1997, LexcellentLicht1991, NguyenMoumni1995, Modelisation1998}.

Let us also note the very important role of internal constraints and of the
reaction ${\bf B}^{ndr}$ to those internal constraints which are responsible
for many properties.

\section{Main results}
\label{mainres}

In order to give a precise formulation of our problem, let us denote by $\Omega$
a bounded, convex set in $\RR^n$ ($n=1,2,3$) with Lipschitz boundary $\Gamma$, by $T$ a positive final time, and
by $Q$ the space-time cylinder $\Omega\times (0,T)$. Let $(\Gamma_0,\Gamma_1)$ be a partition of $\partial\Omega$ into two measurable
sets such that both $\Gamma_0$ and $\Gamma_1$ have positive surface measure. Finally, denote by
$\Sigma:=\partial\Omega\times[0,T]$, $\Sigma_j:=\Gamma_j\times [0,T]$
($j=0,1$) and introduce the Hilbert
triplet $(V,H,V')$ where
\begin{equation}\label{spazi}
H:=L^2(\Omega)\quad\hbox{and}\quad V:=W^{1,2}(\Omega)
\end{equation}
and identify, as usual, $H$ (which stands either for the space $L^2(\Omega)$ or for $(L^2(\Omega))^2$ or for
$(L^2(\Omega))^3$) with its dual space $H'$, so that $V\hookrightarrow H\hookrightarrow V'$ with dense
and continuous embeddings. Moreover, we denote by $\Vert\cdot\Vert_X$ the norm in some space $X$ and by
$\duav{\cdot,\cdot}$ the duality pairing between $V$ and $V'$ and by $(\cdot,\cdot)$ the scalar product in $H$.

Set, for simplicity of notation and without any loss of generality
$$c=k=l_a/\teta_0=\bar C=\lambda=\upsilon=1\,.$$

Then, in order to write the variational formulation of our problem (\ref{pde1}--\ref{pde5}),
we need to generalize the relationship ${\bf B}^{ndr}\in \partial I_{C}(\betab)$ stated
in \eqref{pde4} (cf.~also \cite{bcgg} for similar generalizations). Hence, we need to introduce
the following ingredients
\begin{align}\no
&j:\,\RR^3\to [0,+\infty]\hbox{ a proper, convex, lower semicontinuous function}\\
\label{jpiccola}
&\qquad\hbox{ such that }j(0)=0\quad\hbox{and its subdifferential}\\
\label{alfa}
&\alpha=\partial j:\,\RR^3\to 2^{\RR^3}\\
\label{tau}
&\tau\in W^{1,\infty}(\RR)\,.
\end{align}
Moreover, we consider the associate functionals
\begin{align}\label{jh1}
&J_H(\vb)=\io j(\vb(x))dx\quad\hbox{if}\quad\vb\in H\quad\hbox{and}\quad j(\vb)\in L^1(\Omega)\\
\label{jh2}
&J_H(\vb)=+\infty\quad\hbox{if}\quad\vb\in H\quad\hbox{and}\quad j(\vb)\not\in L^1(\Omega)\\
\label{jv}
&J_V(\vb)=J_H(\vb)\quad\hbox{if}\quad\vb\in V^3
\end{align}
with their subdifferentials
(cf.~\cite[Chap.~II, p.~52]{barbu})
\begin{equation}\label{subvp}
\partial_{V,V'}J_V:\,V^3\to 2^{(V')^3}
\end{equation}
and (cf.~\cite[Ex.~2.1.4, p.~21]{brezis})
\begin{equation}\label{subh}
\partial_H J_H:\,H\to 2^{H}.
\end{equation}
Denote by $D(\partial_{V,V'} J_V):=\{\vb\in V^3\,:\,\partial_{V,V'} J_V(\vb)\neq\emptyset\}$
the domain of $\partial_{V,V'}J_V$. Then, for $\chib,\,\xib\in H$, we have
(see, e.g., \cite[Ex.~2.1.3, p.~52]{brezis}) that
\begin{equation}\no
\xib\in\partial_H J_H(\chib)\quad\hbox{if and only if }\xib\in\alpha(\chib)\quad\hbox{a.e. in }\Omega
\end{equation}
and, thanks to \eqref{jv} and to the definitions of $\partial_{V,V'}J_V$ and $\partial_H J_H$, we have
\begin{equation}\label{inclu}
\partial_H J_H(\chib)\subseteq H\cap\partial_{V,V'} J_V(\chib)\quad\forall\chib\in V^3.
\end{equation}
Now we denote by $W$ the following space
\begin{equation}\label{Vspace}
W:=\{\vb\in V^3\,:\,\vb={\bf 0}\hbox{ on }\Gamma_0\}
\end{equation}
endowed with the usual norm. In addition, we introduce on $W\times W$ a bilinear symmetric
continuous form $a(\cdot,\cdot)$
defined by
\begin{equation}\no
a(\ub,\vb):=\sum_{i,j=1}^3\io\e_{ij}(\ub)\e_{ij}(\vb).
\end{equation}
Note here that (since $\Gamma_0$ has positive measure), thanks to
Korn's inequality (cf., e.g., \cite{ciarlet}, \cite[p.~110]{DL}),
 there exists a positive constant $c$ such that
\begin{equation}\label{aequiv}
a(\vb,\vb)\geq c\Vert\vb\Vert_W^2\quad\forall\vb\in W.
\end{equation}
Next, in order to rewrite the problem (\ref{pde1}--\ref{pde5}) in an abstract framework,
let us introduce the operators
\begin{align}\label{laplneu}
&{\cal B}:V\to V',\quad\duav{{\cal B}u, v}=\io\nabla u\cdot\nabla v \quad u,v\in V\\
\label{opea}
&{\cal A}:W\to W',\quad\duaw{{\cal A}\ub,\vb}=a(\ub,\vb)\quad \ub,\vb\in W\\
\label{opediv}
&{\cal H}:H\to W',\quad \duaw{{\cal H}u,\vb}=\io u\,\dive\vb\quad u\in H,\vb\in W.
\end{align}
Moreover, we make the following assumptions on the data
\begin{align}\label{datou}
&\ub_0\in W,\quad\teta_0\in L^1(\Omega),\\
\label{datoteta}
&\teta^0>0\quad\hbox{a.e. in }\Omega,\quad w_0:=\log\teta^0\in W^{2,2}(\Omega),\\
\label{datochi}
&\betab_0=(\beta_1^0,\beta_2^0,\beta_3^0)\in D(\partial_{V,V'}J_V),\\
\label{sorg1}
&{\bf f}\in W^{1,1}(0,T;H),\quad {\bf g}\in W^{1,1}(0,T;(L^2(\Gamma_1))^3),\\
\label{sorg2}
&R\in H^{1,1/2}(Q)\cap L^1(0,T;L^\infty(\Omega)),\quad \Pi\in L^{\infty}(\Sigma)\cap H^{3/2,3/4}(\Sigma).
\end{align}
Then, we introduce the functions ${\cal R}\in\LDVp$ and ${\cal F}\in W^{1,1}(0,T, W')$ such that
\begin{align}\label{sorg3}
&\duav{{\cal R}(t),v}=\io R(t)v+\int_{\partial\Omega} \Pi(t) v_{|\partial\Omega}\quad v\in V,
\quad \hbox{for a.e. } t\in [0,T],\\
\label{sorg4}
&\duaw{{\cal F}(t),\vb}=\io {\bf f}(t)\cdot\vb
+\int_{\Gamma_1}{\bf g}(t)\cdot\vb_{|\partial\Omega}\quad\vb\in W,\quad\hbox{for a.e. }t\in [0,T].
\end{align}

Now take the function
\begin{equation}\label{gammateta}
\gamma(r):=\exp(r)\quad\hbox{for }r\in \RR
\end{equation}
and  take $J_V$ as in \eqref{jv}, then we are ready to introduce the variational formulation
of our problem as follows.

\medskip
\noindent
{\sc Problem (P).} Given $\teta_0>0$
find $(\ub,\,w,\,\beta_1,\beta_2,\beta_3)$ and $(\xi_1,\xi_2,\xi_3, \,p)$ with the regularities
\begin{align}\label{reg1}
&\ub\in L^\infty(0,T, W),\quad\dive\ub\in \HUV\\
\label{reg2}
&w\in \HUH\cap\LDV\cap L^\infty(Q)\\
\label{reg3}
&\beta_1,\beta_2, \beta_3\in \HUV\cap \LIV,\quad \betab\in D(\partial_{V,V'} J_V)\quad\hbox{a.e.~in }(0,T)\\
\label{reg4}
&\xi_1,\xi_2,\xi_3\in\LDVp,\quad p\in L^2(Q)
\end{align}
satisfying
\begin{align}\label{p1}
&{\cal A}\ub-{\cal H}(p+(\beta _{1}-\beta _{2})\tau (\gamma(w)))={\cal F}\quad\hbox{in }W'\hbox{ a.e.~in }[0,T]\\
\label{p2}
&\dt(\beta_1+\beta_2+\beta_3)+\dive \ub_t=0\quad\hbox{a.e.~in }Q\\
\label{p3}
&\dt w+\dt(\beta_3)+{\cal B}w={\cal R}\quad\hbox{in }V'\hbox{ a.e.~in }[0,T]\\
\no
&\dt\betab+\begin{pmatrix}{\cal B}\dt{\beta_1}\\{\cal B}\dt{\beta_2}\\{\cal B}\dt{\beta_3}\end{pmatrix}
+\begin{pmatrix}{\cal B}\beta_1\\{\cal B}{\beta_2}\\{\cal B}{\beta_3}\end{pmatrix}+\xib
=\begin{pmatrix}
-\tau (\gamma(w)):\varepsilon(\ub) -p \\
\tau (\gamma(w)):\varepsilon(\ub) -p \\
-(\gamma(w)-\teta_0)-p
\end{pmatrix}\\
\label{p4}
&\qquad\hbox{in }(V')^3\hbox{ a.e.~in }[0,T]
\\
\label{p5}
&\xib=(\xi_1,\xi_2,\xi_3)\in\partial_{V,V'}J_V(\betab)\quad\hbox{a.e.~in }[0,T]
\end{align}
and such that
\begin{align}\label{p6}
&\ub(0)=\ub_0\quad\hbox{a.e.~in }\Omega\\
\label{p7}
&w(0)=w_0\quad\hbox{a.e.~in }\Omega\\
\label{p8}
&\betab(0)=\betab_0\quad\hbox{a.e.~in }\Omega\,.
\end{align}

\beos\label{kappa} Obviously we can take $j=I_{C}$ in
\eqref{jpiccola} with $C$ as in \eqref{defC} (cf.~\cite[Ex.~3,
p.~54]{barbu}) and recover the problem already stated in
(\ref{pde1}--\ref{pde5}) as a particular case of our more general
formulation. Note moreover that we can write down equation
\eqref{p4} only in $(V')^3$ (and by consequence we need to
introduce the notion \eqref{subvp} of subdifferential in $(V')^3$)
because of the $V'$ regularity of $({\cal B}\dot\beta_1,{\cal
B}\dot\beta_2,{\cal B}\dot\beta_3)$ in \eqref{p4}. The difficult
point in the proof of our result will be in fact the passage to
the limit in the two non-smooth nonlinearities in equation
\eqref{p4}. By the contrary, we aim to remark that this viscous
term in \eqref{p4} gives more spatial regularity to $\dt\betab$
which furnish more regularity to $w$ in \eqref{p3} (cf.~the
following Theorem~\ref{exiteo}) and consequently to
$\teta=\gamma(w)$ in \eqref{p4}. This regularity is needed in
order to prove well-posedness for our problem. However, we can
also notice that, from the mechanical viewpoint, since we have
introduced in the model the elastic (non-dissipative) local
 interaction term $-\Delta \betab$ (in \eqref{pde2}), it
seems also reasonable to include in the model the dissipative
local
interaction term $-\Delta\betab_t$. \eddos

We are now ready to state our main result which is the
following global existence and uniqueness theorem.
\bete\label{exiteo}
Let  the assumptions on the data (\ref{jpiccola}--\ref{sorg4})
hold and let $T$ be a positive final time.
Then {\sc Problem (P)} has a unique solution on the whole time interval $[0,T]$.
\ente
Moreover, in the last Section~\ref{prouni}, we will get a proof for the following
continuous dependence result
for {\sc Problem (P)}.
\bete\label{uniteo}
Let $T$ be a positive final time, $(\ub^i_0,\,w_0^i,\,\betab_0^{i})$ ($i=1,2$)
be two sets of initial data satisfying conditions (\ref{datou}--\ref{datochi}),
$({\cal R}^i, \,{\cal F}^i)$ ($i=1,2$) be two data of {\sc Problem (P)} satisfying
assumptions (\ref{sorg3}--\ref{sorg4}) with
$({\bf f}^i,\,{\bf g}^i,\,R^i,\,\Pi^i)$ ($i=1,2$) as in (\ref{sorg1}--\ref{sorg2}).
Let $(\ub^i,\,w^i,\,\beta_1^i,\,\beta_2^i,\,\beta_3^i)$
($i=1,2$) be two solutions of {\sc Problem (P)} corresponding to these data.
Moreover, besides conditions (\ref{jpiccola}--\ref{alfa}),
suppose that the following hypothesis
\begin{equation}\label{ipoalfalip}
\alpha\in C^{0,1}(\RR^3)
\end{equation}
holds. Then, there exists a positive constant $M$, depending on the data of the problem, such that the
following continuous dependence estimate
\begin{align}\no
&\| w^1-w^2\|_{\LItH\cap\LDtV}^2+\|\ub^1-\ub^2\|_{\HUtW}^2+\sum_{j=1}^3\|\beta_j^1-\beta_j^2\|_{\HUtV}^2\\
\no
&\leq M\Big(\|w_0^1-w_0^2\|_H^2+\|{\bf u}_0^1-{\bf u}_0^2\|_W^2+\sum_{j=1}^3\|\beta_j^{01}-\beta_j^{02}\|_V^2\\
\label{CD}
&\quad\qquad+\|{\cal F}^1-{\cal F}^2\|_{L^2(0,t;W')}^2+\|{\cal R}^1-{\cal R}^2\|_{\LDtVp}^2\Big)
\end{align}
holds for any $t\in (0,T)$.
\ente

\beos\label{uni}
Let us note that in this paper we can treat
the difficult coupling between the phase-equations \eqref{p4} in which
appears the temperature $\teta$ ($=\gamma(w)$) and the entropy balance equation \eqref{p3} in which
only the function $\log\teta$ ($=w$) plays some role, using the $L^\infty$-bound on the
$w$-component of solution to {\sc Problem (P)}
(cf.~\eqref{reg2}). Indeed it was just due to the lack of regularity of solutions that in \cite{bcf}
(where there was the same type of coupling without the $\Delta{\dt \beta}$-term in \eqref{p4})
the authors did not obtain uniqueness of solutions (cf.~also \cite[Remark~5.2]{bcf}).

However, let us, finally, observe that the main advantage of taking the entropy balance equation
instead of the internal energy balance equation is that once one has solved the problem in some sense and has found
the temperature $\teta:=\gamma(w)$, it is automatically positive because it stands
in the image of the function $\gamma$
(cf.~\eqref{gammateta}). Indeed in many cases it is difficult to deduce this fact only
from the internal energy balance
equation  (cf., e.g., \cite{cspos} in order to see one example of these difficulties).
Let us note that within the small
 perturbations assumption the entropy balance and the classical heat equation are
 equivalent in mechanical terms
 (cf.~\cite{bcf, bf, bfr}).
\eddos

\section{Proof of Theorem~\ref{exiteo}}
\label{exiproof}

The following section is devoted to the proof of
Theorem~\ref{exiteo}. First we
approximate our {\sc Problem (P)} by a more regular {\sc Problem $(P\ee)$},
then (fixed $\e>0$) we find well-posedness for the approximating
problem using a iterated Banach contraction fixed-point argument and then we perform some
a-priori estimates (independent of $\e$) on its solution, which allow us
to pass to the limit  in {\sc Problem $(P\ee)$} as $\e\searrow 0$, recovering a solution to {\sc Problem (P)}.

\subsection{The approximating problem}\label{appro}

We take a small positive parameter $\e>0$ and approximate $\partial_{V,V'}J_V$ in (\ref{p4}--\ref{p5})
, let us take the Lipschitz continuous Yosida-Moreau approximation $\alpha\ee=\partial j\ee=(j\ee)'$
(cf.~\cite[Prop.~2.11, p.~39]{brezis}) of $\alpha$ and the associated functional $J_{H,\e}(\vb)=\io j\ee(\vb(x))\,dx$,
whose differential ($\partial_H J_{H,\e}$) is the Yosida-Moreau approximation of $\partial_H J_H$ (cf.~\cite[Prop.~2.16,
p.~47]{brezis}). Now we aim to recall some properties of this approximation which will be useful in order to pass to
the limit as $\e\searrow 0$. Note that the proof of the following lemma is a consequence of
\cite[Thm.~3.20, p.~289]{attouch}, \cite[Thm.~3.62, p.~365]{attouch}, and of the Lebesgue theorem of passage to the limit
under the sign of integral.
\bele\label{projep}
If $\partial_{V,V'}J_{V,\e}$ is the Yosida-Moreau approximation of $\partial_{V,V'} J_V$ (cf.~\cite[Thm.~2.2, p.~57]{barbu}),
then the following inclusion
\begin{equation}\label{pj1}
\partial_H J_{H,\e}\subseteq\partial_{V,V'}J_{V,\e}
\end{equation}
holds true. Moreover, the following properties hold true (for $\e\searrow 0$)
\begin{align}\label{l1}
&J_{V,\e}\to\sup_{\e>0}J_{V,\e}\hbox{ in the sense of Mosco (cf.~\cite[Def.~3.17, p.~295]{attouch})},\\
\label{l2}
&\sup_{\e>0}J_{V,\e}=\lim_{\e\searrow 0}\io j\ee=\io j,\quad\hbox{hence}\\
\label{l3}
&J_{V,\e}\to J_V\hbox{ in the sense of Mosco and}\\
\label{l4}
&\forall\, \vb\in V^3,\,\,\forall\,\vb_{\e}\to\vb \hbox{ weakly in }V^3\,\,J_V(\vb)<\liminf J_{V,\e}(\vb_{\e}).
\end{align}
Finally, for $\e\searrow 0$, it holds
\begin{equation}\label{pj2}
(\ub,\partial_{V,V'}J_{V,\e}(\ub))\to\partial_{V,V'}J_V\hbox{ in the graph sense (
cf.~\cite[Def.~3.58, p.~360]{attouch}).}
\end{equation}
\enle

Then, let us call $\gamma\ee$ the following Lipschitz continuous
approximation of the function $\gamma(w)=\exp(w)$, i.e.~the
function
\begin{equation}\label{gamma}
\gamma\ee(r):=\begin{cases}
                    \exp r&\hbox{if }r\leq 1/\e\\
                    (r-1/\e)\exp(1/\e)+\exp(1/\e)&\hbox{if }r\geq 1/\e.
                 \end{cases}
\end{equation}
Moreover let $\delta\ee$ be the inverse function of $\gamma\ee$, i.e.
\begin{equation}\label{delta}
\delta\ee\big(\gamma\ee(r)\big)=r\quad\forall r\in\RR
\end{equation}
and let $\widehat\gamma\ee$ be a primitive of the function $\gamma\ee$, i.e.
\begin{equation}\label{gammaprim}
\widehat\gamma\ee(r)=1+\int_0^r\gamma\ee(s)\,ds\quad\forall r\in \RR.
\end{equation}
Then the following properties of $\gamma\ee$ hold true (cf.~also \cite[Lemma~5.1]{bcf}).
\bele\label{progamma}
There holds
$$\widehat\gamma\ee(r)\geq\gamma\ee(r),\quad r(\delta\ee)'(r)\geq 1\quad\forall r\in\RR.$$
\enle
We are ready now to introduce the approximating {\sc Problem $(P\ee)$} as follows.
\smallskip

\noindent
{\sc Problem $(P\ee)$.}
Given $\teta_0>0$
find $(\ub\ee,\,w\ee,\,\beta_1\ee,\beta_2\ee,\beta_3\ee,p\ee)$   with the regularities
\begin{align}\label{reg1e}
&\ub\ee\in L^\infty(0,T, W),\quad\dive\ub\ee\in\HUV, \quad p\ee\in\LDH\\
\label{reg2e}
&w\ee\in \HUH\cap\LDV\cap L^\infty(Q)\\
\label{reg3e}
&\beta_1\ee,\beta_2\ee, \beta_3\ee\in \HUV
\end{align}
satisfying
\begin{align}\label{p1e}
&{\cal A}\ub\ee-{\cal H}(p\ee+(\beta _{1}\ee-\beta _{2}\ee)\tau (\gamma\ee(w\ee)))={\cal F}\quad\hbox{in }W'\hbox{ a.e.~in }[0,T]\\
\label{p2e}
&\dt(\beta_1\ee+\beta_2\ee+\beta_3\ee)+\dive (\ub\ee)_t=\varepsilon p\ee\quad\hbox{a.e.~in }Q\\
\label{p3e}
&\dt w\ee+\dt\beta_3\ee+{\cal B}w\ee={\cal R}\quad\hbox{in }V'\hbox{ a.e.~in }[0,T]\\
\no
&\dt\betab\ee+\begin{pmatrix}{\cal B}\dt{\beta_1\ee}\\{\cal B}\dt{\beta_2\ee}\\{\cal B}\dt{\beta_3\ee}\end{pmatrix}
+\begin{pmatrix}{\cal B}\beta_1\ee\\{\cal B}{\beta_2\ee}\\{\cal B}{\beta_3\ee}\end{pmatrix}+\partial_H J_{H,\e}(\betab\ee)
=\begin{pmatrix}
-\tau (\gamma\ee(w\ee)):\varepsilon(\ub\ee) -p\ee \\
\tau (\gamma\ee(w\ee)):\varepsilon(\ub\ee) -p\ee \\
-(\gamma\ee(w\ee)-\teta_0)-p\ee
\end{pmatrix}
\\
\label{p4e}
&\qquad\hbox{in }(V')^3\hbox{ a.e.~in }[0,T]
\end{align}
and such that
\begin{align}\label{p6e}
&\ub\ee(0)=\ub_0\quad\hbox{a.e.~in }\Omega\\
\label{p7e}
&w\ee(0)=w_0\quad\hbox{a.e.~in }\Omega\\
\label{p8e}
&\betab\ee(0)=\betab_0\quad\hbox{a.e.~in }\Omega\,.
\end{align}

Concerning this approximating problem, we prove hereafter
the following well-posedness result.
\bete\label{teoep}
Let the assumptions (\ref{jpiccola}--\ref{sorg4}) hold true. Let $T$ be a positive final
time and $\e>0$. Then the {\sc Problem $(P\ee)$} has a unique solution in $[0,T]$.
\ente
\prova
Here we are going first to prove local existence (and uniqueness) in a finite time interval $[0,\bar{t}]$ for some
$\bar{t}\in [0,T]$, then we will extend the solution to the whole interval $[0,T]$ proving
global existence (and uniqueness) of solution to {\sc Problem $(P\ee)$}. Hence, let us take $\bar{t}\in [0,T]$
(we will choose it later) and denote by ${\cal X}:=(W^{1,2}(0,\bar{t};V))^3$.
Fix for the moment $(\betapb\ee,\betasb\ee,\betatb\ee)\in {\cal X}$ in the equations (\ref{p1e}--\ref{p3e}), then,
by well-known results (cf.~also (\ref{sorg1}--\ref{sorg2})), we find  a unique
$w\ee={\cal T}_2(\betapb\ee,\betasb\ee,\betatb\ee)\in \HUH\cap\CZV\cap L^2(0,T;W^{2,2}(\Omega))$  solution of \eqref{p3e}. By
\cite[Thm.~6.2, p.~168]{DL}, it is possible to find a unique $\ub\ee={\cal T}_1(\betapb\ee,\betasb,\betatb\ee,w\ee)\in L^\infty(0,T;W)$
solution of (\ref{p1e}--\ref{p2e}) such that $\dive (\ub\ee)_t\in L^2(0,T;H)$.

Moreover, if we take these values of $\ub\ee$ and $w\ee$ in \eqref{p4e},
we can find a solution $(\beta_1\ee,\beta_2\ee,\beta_3\ee)\in {\cal X}$
(depending on $\ub\ee=:{\cal T}_1(\betapb\ee,\betasb\ee,\betatb\ee,w\ee)$ and $w\ee=:{\cal T}_2(\betapb\ee,\betasb\ee,\betatb\ee)$)
of the equation \eqref{p4e} again by standard results.

In this way, we have defined an operator ${\cal T}:\,{\cal
X}\to{\cal X}$ such that $(\beta_1\ee,\beta_2\ee,\beta_3\ee)=:{\cal
T}(\betapb\ee,\betasb\ee,\betatb\ee)$. What we have to do now is to prove that
${\cal T}$ is a contraction mapping on ${\cal X}$ for a
sufficiently small $\bar{t}\in [0,T]$ and moreover, repeating the
procedure step by step in time (this is possible thanks to the regularities
properties of the solution listed above), we can prove well-posedness for
the {\sc Problem $(P\ee)$} on the whole time interval $[0,T]$ and
conclude the proof of Theorem~\ref{teoep}. In order to prove that
${\cal T}$ is contractive, let us proceed by steps and forget of
the apices $\e$.
\smallskip

\noindent
{\sl First step.} Let $(\betapb^i,\betasb^i,\betatb^i)\in{\cal X}$,
$\ub^i={\cal T}_1(\betapb^i,\betasb^i,\betatb^i)$, $w^i={\cal T}_2(\betapb^i,\betasb^i,\betatb^i)$,
and $(\betapb^i,\betasb^i,\betatb^i)={\cal T}(\betapb^i,\betasb^i,\betatb^i)$ ($i=1,2,3$). Then,
writing two times \eqref{p3e} with $\betatb^i$ ($i=1,2$) instead of $\dt(\beta_3)$, making the difference,
testing the resulting  equation with $(w^1-w^2)_t$, and integrating on $(0,t)$ with $t\in [0,T]$,
we get the following inequality
\begin{equation}\label{contr1}
\Vert (w^1-w^2)_t\Vert_{\LDtH}^2+\Vert (w^1-w^2)(t)\Vert_{V}^2
\leq C_0\Vert(\beta_3^1)_t-(\beta_2^2)_t\Vert_{\LDtH}^2+\mezzo\|w_1-w_2\|_{\LDtH}^2,
\end{equation}
for some positive constant $C_0$ independent of $t$. Hence, we get, for all $t\in [0,T]$,
\begin{equation}\label{contr1bis}
\Vert (w^1-w^2)_t\Vert_{L^2(0,t;H)}^2+\frac12\Vert w^1-w^2\Vert_{C^0([0,t];V)}^2
\leq C_1\Vert(\beta_3^1)_t-(\beta_2^2)_t\Vert_{L^2(0,t;H)}^2
\end{equation}
being $C_1:=C_0 {\rm e}^{T/2}$.
\smallskip

\noindent
{\sl Second step.} Let us take $\teta^i=\gamma\ee(w^i)$ and
$\e p^i:=\dt(\beta_1^i+\beta_2^i+\beta_3^i)+\dive (\ub^i)_t$ ($i=1,2$) and write \eqref{p1e} with $\beta_1^i$ and $\beta_2^i$,
make the difference, test the resulting  equation with $(\ub^1-\ub^2)_t$,
integrate on $(0,t)$ with $t\in [0,T]$, and use equation \eqref{p2e}, getting the following inequality
\begin{align}\no
\displaystyle\Vert \ub^1(t)&-\ub^2(t)\Vert_W^2+\varepsilon\Vert p^1-p^2\Vert_{\LDtH}^2\\
\no
&\leq\int_Q\sum_{j=1}^3(\beta_j^1-\beta_j^2)_t(p^1-p^2)\\
\no
&\displaystyle
\quad-\int_Q\left(\tau(\teta^1)\left((\beta_1^1-\beta_1^2)-(\beta_2^1-\beta_2^2)\right)
\right)\left(\varepsilon (p^1-p^2)-\sum_{j=1}^3(\beta_j^1-\beta_j^2)_t\right)\\
\no
&\quad+\int_Q\left(\left(\tau(\teta^1)-\tau(\teta^2)\right)(\beta_1^2-\beta_2^2)\right)\left(\varepsilon (p^1-p^2)
-\sum_{j=1}^3(\beta_j^1-\beta_j^2)_t\right)\\
\no
&\displaystyle\leq \int_Q\sum_{j=1}^3(\beta_j^1-\beta_j^2)_t(p^1-p^2)+\frac{\e}{2}\Vert p^1-p^2\Vert_{\LDtH}^2\\
\no
&\quad+
\frac14\sum_{j=1}^3\Vert(\beta_j^1)_t-(\beta_j^2)_t\Vert_{\LDtH}^2+C_2t^2\e\sum_{j=1}^3\Vert(\beta_j^1)_t-(\beta_j^2)_t\Vert_{\LDtH}^2\\
\label{contr2}
&\displaystyle\quad+C_3 t\exp(2/\e)\Vert w^1- w^2\Vert_{\CZtV}^2\,.
\end{align}

\smallskip
\noindent
{\sl Third step.} Write equation \eqref{p4e} for $\teta^i$ and $\ub^i$,
make the difference between the two equations written for $i=1$ and $i=2$, and test the resulting vectorial equation
by the vector $\big((\beta_1^1)_t-(\beta_1^2)_t, (\beta_2^1)_t-(\beta_2^2)_t, (\beta_3^1)_t-(\beta_3^2)_t\big).$ Summing up the two lines and
integrating on $(0,t)$ with $t\in [0,T]$, we have (exploiting the Lipschitz continuity of $\partial_H J_{H,\e}$
(cf.~\cite[Prop.~2.6, p.~28]{brezis}))
\begin{align}\no
\sum_{j=1}^3\Vert(\beta_j^1)_t&-(\beta_j^2)_t\Vert_{\LDtV}^2
+\sum_{j=1}^3\Vert\nabla(\beta_j^1-\beta_j^2)(t)\Vert_H^2\\
\no
&\leq -\int_Q\sum_{j=1}^3(\beta_j^1-\beta_j^2)_t(p^1-p^2)\\
\no
&\quad-\int_Q\left((\tau(\teta^1)-\tau(\teta^2))\dive\ub^1+(\dive \ub^1-\dive \ub^2)\tau(\teta^2)\right)(\beta_1^1-\beta_1^2)_t\\
\no
&\quad+\int_Q\left((\tau(\teta^1)-\tau(\teta^2))\dive\ub^1+(\dive \ub^1-\dive \ub^2)\tau(\teta^2)\right)(\beta_2^1-\beta_2^2)_t\\
\no
&\quad-\int_Q (\teta^1-\teta^2)(\beta_3^1-\beta_3^2)_t+\frac14\sum_{j=1}^3\Vert(\beta_j^1)_t-(\beta_j^2)_t\Vert_{\LDtH}^2\\
\no
&\quad+\frac{1}{\e^2}\sum_{j=1}^3\Vert\beta_j^1-\beta_j^2\Vert_{\LDtH}^2\,.
\end{align}\
Moreover, using the definition \eqref{progamma} of $\gamma\ee$ and the assumption \eqref{tau} on $\tau$, we get the following inequality
\begin{align}\no
\mezzo\sum_{j=1}^3\Vert(\beta_j^1)_t-(\beta_j^2)_t\Vert_{\LDtV}^2
&+\sum_{j=1}^3\Vert\nabla(\beta_j^1-\beta_j^2)(t)\Vert_H^2\\
\no
&\leq -\int_Q\sum_{j=1}^3(\beta_j^1-\beta_j^2)_t(p^1-p^2)\\
\no
&
\quad+C_4 t\exp(2/\e)\Vert w^1- w^2\Vert_{\CZtV}^2\\
\no
&\quad+C_5\|\ub^1-\ub^2\|_{L^2(0,T;W)}^2+t\exp(2/\e)\Vert w^1- w^2\Vert_{\CZtH}^2\\
\label{contr3}
&\quad+\frac{t^2}{\e^2}\sum_{j=1}^3\Vert(\beta_j^1)_t-(\beta_j^2)_t\Vert_{\LDtV}^2\,.
\end{align}

\smallskip
\noindent
{\sl Fourth step.} Summing up the two inequalities \eqref{contr2} and \eqref{contr3}, two integrals cancels out,
 and using \eqref{contr1}, we get
\begin{align}\no
\sum_{j=1}^3\Vert(\beta_j^1)_t-(\beta_j^2)_t\Vert_{\LDtV}^2&+\sum_{j=1}^3\Vert\nabla(\beta_j^1
-\beta_j^2)(t)\Vert_H^2+\Vert \ub^1-\ub^2\Vert_{C^0([0,t];W)}^2\\
\label{contr4}
&\leq C_\e\big(t+t^2\big)\sum_{j=1}^2\Vert(\beta_j^1)_t-(\beta_j^2)_t\Vert_{\LDtV}^2\,
\end{align}
where $C_\e$ does not depend on $t$. Hence, choosing $t$ sufficiently small (this is our ${\bar t}$),
we recover the contractive property of ${\cal T}$. Moreover,
applying the Banach fixed point theorem to ${\cal T}$, we get a
unique solution for the {\sc Problem $(P\ee)$} on the time interval
$[0,\bar t]$. Due to this estimate it is easy to prove that there exists $m\in \NN$
such that ${\cal T}^m$ is a contraction on $X$. Hence we
have a unique solution on the whole time interval $[0,T]$. This
concludes the proof of Theorem~\ref{teoep}.\fin

\subsection{A priori estimates}\label{esti}

In this subsection we perform a-priori estimates on {\sc Problem $(P\ee)$} uniformly in $\e$ which will lead us
pass to the limit as $\e\searrow 0$ and recover a solution of {\sc Problem (P)}.
We denote by $c$ all the positive constants (which may also differ from line to line) independent
of $\e$ and depending on the data of the problem. For simplicity, we omit the subscript $\e$ when it is not necessary.

\noindent
{\sl First a-priori estimate.} Test \eqref{p1e} by $\ub_t$, \eqref{p2e} by $p$, \eqref{p3e} by $\gamma\ee(w)+w$,
\eqref{p4e} by $\betab_t$, sum up the resulting equations and integrate over $(0,t)$ ($t\in [0,T]$). The result is
\begin{eqnarray}\no
&\displaystyle\io\left(\widehat\gamma\ee(w(t))+\frac12 w^2(t)\right)+\int_Q\nabla\delta\ee(\teta\ee)\nabla\teta\ee+
\int_Q|\nabla w|^2+\itt\|\betab_t\|_V^2
\\
\no
&\displaystyle+\frac12\io|\nabla\betab(t)|^2+J_{H,\e}(\betab(t))+\frac12\|\ub(t)\|_W^2+\e\|p\|_{\LDtH}^2=\io\left(\widehat\gamma\ee(w_0)
+\frac12w_0^2\right)
\\
\no
&\displaystyle+\frac12\|\nabla\betab_0\|_H^2+\frac12\|\ub_0\|_W^2+\itt\,\duaw{{\cal F},\ub_t}+\itt\duav{{\cal R},(\gamma\ee(w)+w)}
\\
\no
&\displaystyle-\int_Q\dt\beta_3 (w-\teta_0)+\int_Q\tau(\gamma\ee(w))\left(\dive\ub\left(\dt\beta_1-\dt\beta_2\right)+(\beta_1-\beta_2)\dive\ub_t\right)\,.
\end{eqnarray}
Now, following the line of \cite[(5.5)--(5.7), p.~1583]{bcf}, we can deal with the source term ${\cal R}$ recalling
\eqref{sorg3} and using a well-known compactness inequality (cf.~\cite[Theorem~16.4]{lima}) in this way
\begin{align}\label{erre}
\int_Q R\gamma\ee(w)&\leq\itt\Vert R(s)\Vert_{L^{\infty}(\Omega)}\Vert\gamma\ee(w(s))\Vert_{L^1(\Omega)}\,ds\\
\no
\itt\int_{\partial\Omega} \Pi\gamma\ee(w)&\leq c\Vert \Pi\Vert_{L^{\infty}(\Sigma)}\Vert (\gamma\ee(w))^{1/2}\Vert_{L^2(\Sigma)}\\
\label{acca}
&\leq \frac12\Vert\nabla(\gamma\ee(w))^{1/2}\Vert_{(\LDtH)^3}^2+c\Vert(\gamma\ee(w))^{1/2}\Vert_{\LDtH}^2\,.
\end{align}
Moreover, using assumptions \eqref{sorg1} and \eqref{sorg4}, we get (integrating
by parts in time)
\begin{align}\no
\itt \,\duaw{{\cal F},\ub_t}&\,=-\itt \,\duaw{{\cal F}_t,\ub}+\, \duaw{{\cal F}(t),\ub(t)}-\, \duaw{{\cal F}(0),\ub(0)}\\
\label{accabis}
&\leq c+\frac14\|\ub(t)\|_W^2+\itt\|{\cal F}_t\|_{W'}\|\ub\|_W\,.
\end{align}
Now, collecting estimates (\ref{erre}--\ref{accabis}), using Lemma~\ref{progamma} and
\eqref{p2e} in order to estimate the term containing $\dive\ub_t$ and employing assumptions
\eqref{tau} on $\tau$ and (\ref{datou}--\ref{sorg4}) on the data,  we get the following inequality
\begin{eqnarray}\no
&\displaystyle\io\left(\gamma\ee(w(t))+\frac12 w^2(t)\right)+\int_Q\frac{|\nabla\gamma\ee(w)|^2}{\gamma\ee(w)}
+\int_Q|\nabla w|^2+\frac14\itt\|\betab_t\|_V^2+\frac12\|\betab(t)\|_V^2
\\
\no
&\displaystyle+J_{H,\e}(\betab(t))+\frac14\|\ub(t)\|_W^2+\e\|p\|_{\LDtH}^2\leq c+\int_Q w^2
+\itt\|\ub\|_W^2+c\itt\|\betab\|_H^2\\
\no
&\displaystyle+\frac{\e^2}{4}\itt\|p\|_H^2+c\Vert(\gamma\ee(w))^{1/2}\Vert_{\LDtH}^2+\itt\|{\cal F}_t\|_{W'}\|\ub\|_W
\end{eqnarray}
which - via Gronwall lemma - and choosing $\e$ small, leads to the following first estimate
\begin{eqnarray}\no
&\displaystyle\Vert\ub(t)\Vert_W^2+\Vert(\gamma\ee(w))^{1/2}(t)\Vert_H^2+
\itt\Vert\nabla(\gamma\ee(w))^{1/2}\Vert_{H}^2+\|w(t)\|_{H}^2+\|w\|_{\LDtV}^2\\
\label{s1}
&\displaystyle+\e\|p\|_{\LDtH}^2+\Vert\dt\betab\Vert_{\LDtV}^2+\Vert\nabla\betab(t)\Vert_H^2+J_{H,\e}(\betab\ee(t))\leq c\,.
\end{eqnarray}
Then, applying the standard regularity results for linear parabolic equations to \eqref{p3e}, we get
\begin{equation}\label{s2}
\Vert w\Vert_{\HUH\cap \LIV\cap L^2(0,T;W^{2,2}(\Omega))}\leq c\,.
\end{equation}
Now, thanks to assumptions \eqref{datoteta} and \eqref{sorg2}, we may apply \cite[Lemma~4.3]{frm3as}
to equation \eqref{p2e}, getting the further bound
\begin{equation}\label{s3}
\Vert w\Vert_{L^\infty(Q)}\leq c\,.
\end{equation}
Thanks to \eqref{gamma} we can immediately recover
\begin{equation}\label{s4}
\Vert\gamma\ee(w)\Vert_{L^{\infty}(Q)}\leq c
\end{equation}
and, due to Lemma~\ref{progamma} and (\ref{s1}--\ref{s4}), we have
\begin{align}\no
\int_Q\vert\nabla\gamma\ee(w)\vert^2\leq c\int_Q\frac{\vert\nabla\gamma\ee(w)\vert^2}{\gamma\ee(w)}&\leq
c\int_Q\nabla\delta\ee(\gamma\ee(w))\nabla\gamma\ee(w)\\
\no
&\leq c\int_Q\vert\nabla w\vert^2\leq c
\end{align}
and
\begin{align}\no
\int_Q\vert(\gamma\ee(w))_t\vert^2\leq c\int_Q\Big\vert\frac{(\gamma\ee(w))_t}{\gamma\ee(w)}\Big\vert^2
&\leq c\int_Q\vert(\delta\ee)'(\gamma\ee(w))(\gamma\ee(w))_t\vert^2\\
\no
&\leq c\int_Q\vert w_t\vert^2\leq c.
\end{align}
Hence, from these two inequalities  it follows that
\begin{equation}\label{s5}
\Vert\gamma\ee(w)\Vert_{\HUH\cap\LDV\cap L^{\infty}(Q)}\leq c\,.
\end{equation}

\noindent
{\sl Second a-priori estimate.} In order to pass to the limit (as $\e\searrow 0$) in \eqref{p4e},
we need to pass to the limit in $p\ee$. We will use the following Lions'
lemma which is stated in this form, e.g.,  in \cite[Rem.~1.1, p.~17]{templast} (its proof is due to
\cite[Note~(27), p.~320]{ms} in case of a ${\cal C}^1$ class domain $\Omega$ and to \cite{necas}
when $\Omega$ is only Lipschitz). For further comments on this topic the reader can refer
to \cite[Remark~4.1]{frm3as}.
\bele\label{lelions}
Let $\Omega$ be a bounded and Lipschitz set in $\RR^3$ and let $m$ be a continuous seminorm on $H$ and
a norm on the constants. Then there exists a positive constant $c(\Omega)$ (depending only on $\Omega$)
such that the following inequality
\begin{equation}\label{poinc}
\Vert u\Vert_{H}\leq c(\Omega)\{m(u)+\Vert\nabla u\Vert_{V'}\}
\end{equation}
holds for all $u\in H$ with $\nabla u\in V'$.
\enle
We want to apply this result in order to find the uniform (in $\e$) bound on $p\ee$. First of all
let us note that from comparison in \eqref{p1e}, using also the bound \eqref{s1} on $\ub$
with the assumption \eqref{sorg1} on ${\cal F}$, we immediately deduce that
\begin{equation}\label{stimaterza}
\Vert \nabla p\ee\Vert_{L^2(0,T;V')}\leq c\,.
\end{equation}
Moreover, always by comparison in \eqref{p1e}, we have that
\begin{equation}\label{stimaquarta}
\Big\vert\int_Q p\ee\dive \vb\Big\vert\leq c\quad\forall \vb\in W.
\end{equation}
Following the idea of \cite{frm3as}, we can choose $\vb_{\star}\in W$ such that
\begin{equation}\label{ipodom}
\io\dive(\vb_{\star})\,dx=\int_{\partial\Omega} \vb_{\star}\cdot {\bf n}\,ds\neq 0.
\end{equation}
Note that, since $\Omega$ is regular (it suffices for $\Omega$ to
be a Lipschitz domain), we can always find a $\vb_{\star}\in W$
such that \eqref{ipodom} is satisfied, because, if we take
$B_{\e}(x)$ the ball in $\RR^3$ centered in $x\in \Gamma_1$ with
radius $\e$ such that $B_{\e}(x)\cap\Gamma_0=\emptyset$ and
consider the parametrization of $\Gamma_1\cap B_{\e}(x)$ through
the Lipschitz function
$(x_1,x_2)\mapsto(x_1,x_2,\varphi(x_1,x_2))$,
 then the normal unit vector associated is
$${\bf n}=\frac{(\partial_{x_1}\varphi,-\partial_{x_2}\varphi,1)}{\sqrt{1+\vert\nabla\varphi\vert^2}}.$$
Then, if we take $\vb_{\star}=(0,\,0,\,\zeta)$ with
\begin{equation}\no
\zeta(y)=\begin{cases}\exp{\Big(-\frac{1}{1-\frac{\vert x-y\vert^2}{\e^2}}\Big)}&\hbox{if}\quad\vert x-y\vert\leq\e\\
0&\hbox{otherwise,}
\end{cases}
\end{equation}
then $\vb_{\star}\in W$ and moreover we can show that \eqref{ipodom} holds because
$$\frac{1}{\sqrt{1+\vert\nabla\varphi\vert^2}}\geq\frac{1}{\sqrt{1+L^2}},$$
where $L$ is the Lipschitz constant of $\varphi$, and hence
\begin{equation}\no
\io\dive{\vb_{\star}}=\int_{\Gamma_1}\vb_{\star}\cdot{\bf n}\,ds=\int_{\Gamma_1\cap B_{\e}(x)}
\frac{\zeta}{\sqrt{1+\vert\nabla\varphi\vert^2}}\,ds\neq 0\,.
\end{equation}

Take now $m$ in Lemma~\ref{lelions} as
$$m(v)=\Big\vert\io v\dive \vb_{\star}\Big\vert.$$
Then, $m(v)$ is a seminorm on $H$ and a norm on the constants  because of \eqref{ipodom}.
Hence, we can apply Lemma~\ref{lelions} to $p\ee$ with the choices  done above and,
thanks to (\ref{stimaterza}--\ref{stimaquarta}), we get the bound
\begin{equation}\label{s6}
\Vert p\ee\Vert_{L^2(Q)}\leq c\,.
\end{equation}

Finally, by comparison in \eqref{p4e} and using the estimates (\ref{s1}) and (\ref{s5}--\ref{s6}),
we deduce that also $\partial_H J_{H,\e}(\betab)$ is bounded in $\LDVp$.
Then, testing  \eqref{p4e} with ${\cal B}\betab\ee$ and then by using again
(\ref{s1}) and (\ref{s5}--\ref{s6}) and the monotonicity properties of $\alpha\ee$, we get also
\begin{equation}\label{s8}
\Vert\beta_j\ee\Vert_{L^\infty(0,T;W^{2,2}(\Omega))}\leq c\quad (j=1,2,3).
\end{equation}
Now it remains only to pass to the limit in (\ref{p1e}--\ref{p4e}) as $\e\searrow 0$.
This will be the aim of the next subsection.

\subsection{Passage to the limit and uniqueness}\label{conclu}

As we have just mentioned, we want to conclude the proof of Theorem~\ref{exiteo} passing to the limit in the
well-posed (cf.~Subsection~\ref{appro}) {\sc Problem $(P\ee)$} as $\e\searrow 0$ using the previous uniform
(in $\e$) estimates on its solution (cf.~Subsection~\ref{esti}) and exploiting some compactness-monotonicity argument.
Let us list before the weak or weak-star convergence coming directly from the previous estimates and well-known
weak-compactness results. Note that the following convergences hold only up to a subsequence of $\e\searrow 0$
(let us say $\e_k\searrow 0$). We denote it again with $\e$ only for simplicity of notation. From the estimates
(\ref{s1}--\ref{s8}) and the property \eqref{pj1} of $\partial_H J_{H,\e}$, we deduce that
\begin{align}\label{c1}
\ub\ee\to \ub\quad&\hbox{weakly star in }\LIW\\
\label{c2}
\dive \ub\ee_t\to \dive \ub_t\quad&\hbox{weakly in }\LDH\\
\no
w\ee\to w\quad&\hbox{weakly in }W^{1,2}(0,T;H)\cap L^2(0,T; W^{2,2}(\Omega))\\
\label{c3}
&\hbox{and weakly star in }L^\infty(Q)\cap \LIV\\
\label{c4}
\gamma\ee(w\ee)\to\teta\quad&\hbox{weakly star in }\HUH\cap\LIV\cap L^\infty(Q)\\
\no
\beta_j\ee\to\beta_j\quad&\hbox{weakly star in }W^{1,2}(0,T;V)\cap L^\infty(0,T;W^{2,2}(\Omega))\\
\label{c5}
&(j=1,2,3)\\
\label{c6}
p\ee\to p\quad&\hbox{weakly in }L^2(Q)\\
\label{c7}
\xi_j\ee\to \xi_j\quad&\hbox{weakly in }L^2(0,T;V')\quad (j=1,2,3)\,.
\end{align}
Moreover, employing \cite[Cor.~5, p.~86]{simon}, we get also
\begin{align}\no
w\ee\to w\quad&\hbox{strongly in }\LDV\cap\CZH\\
\label{c8}
&\hbox{and hence  a.e.~in }Q\\
\label{c9}
\teta\ee\to\teta\quad&\hbox{strongly in }\CZH\\
\label{c10}
\beta_j\ee\to\beta_j\quad&\hbox{strongly in }\CZV\quad (j=1,2,3)\,.
\end{align}
Note that (\ref{c8}--\ref{c9}) imply immediately the convergence
\begin{equation}\no
\gamma\ee(w\ee)\to\teta=\gamma(w)\quad\hbox{and }\tau(\gamma\ee(w\ee))\to\tau(\teta)\quad\hbox{a.e.~in }Q\,.
\end{equation}
Moreover, the two convergences \eqref{c7} and \eqref{c10} along
with the property \eqref{pj2} and \cite[Thm.~3.66,
p.~373]{attouch} give immediately the identification of the
maximal monotone graph $\partial_{V,V'}J_V$, i.e.
\begin{equation}\no
\xib\in\partial_{V,V'}J_V(\betab)\quad\hbox{in }(V')^3\hbox{ and
a.e. in }[0,T]
\end{equation}
with $\xib=(\xi_1,\xi_2,\xi_3)$ and $\xi_j$ $(j=1,2,3)$ are the weak limits defined in \eqref{c7}.
All these convergences with the identifications made above make us able to pass to the
limit (as $\e\searrow 0$ or at least for a subsequence of it) in {\sc Problem $(P\ee)$}
finding a solution to {\sc Problem (P)} and concluding in this way the proof of Theorem~\ref{exiteo}.
Note that the convergences hold for all subsequences $\e_k$ of $\e$
tending to 0 because of uniqueness of solutions. Indeed we may prove it in this way.

Consider two solutions of {\sc Problem (P)}  $(\ub^i,\,w^i,\,\beta_1^i,\,\beta_2^i,\beta_3^i, p^i)$
($i=1,2,3$) corresponding to the same data. Moreover let us take the {\sl mass balance} equation
\eqref{p2} in the following integrated form
\begin{equation}\label{mi}
\beta_1+\beta_2+\beta_3+\dive \ub=\beta_1(0)+\beta_2(0)+\beta_3(0)+\dive\ub(0)\quad \hbox{a.e.~in }Q\,.
\end{equation}
Then, integrate equation \eqref{p3} over $(0,t)$ (let us call it $1*\eqref{p3}$ with a little abuse of
notation) and write down
two times equations (\ref{p1}--\ref{p2}), $1*\eqref{p3}$, \eqref{p4} with
$(\ub^i,\,w^i,\,\beta_1^i,\,\beta_2^i,\,\beta_3^i, \, p^i)$, make the difference
between the two equations $1*\eqref{p3}$, and
test the result with $(w^1-w^2)$.  Make the difference between the two equations \eqref{p1},
 test the result with $(\ub^1-\ub^2)$. Make the difference between the two equations \eqref{p4},
written for $i=1$ and $i=2$, and test the resulting vectorial equation
by the vector $\left((\beta_1^1)-(\beta_1^2), (\beta_2^1)-(\beta_2^2), (\beta_3^1)-(\beta_3^2)\right).$
 Finally, summing up the three resulting equations, integrating
 over $(0,t)$, with $t\in [0,T]$, exploiting the monotonicity of $\partial_{V,V'}J_V$, using equation \eqref{mi}
in order to get rid of the $p$-terms,
and using the fact that $\gamma$, defined in \eqref{gammateta}, is a locally Lipschitz continuous function,
$\teta^i=\gamma(w^i)$ ($i=1,2$), and $w^i$  are bounded in $L^\infty(Q)$ (cf.~\eqref{reg2}),
we get the following inequality
\begin{eqnarray}\no
&\displaystyle \mezzo\Vert (\ub^1-\ub^2)\Vert_{L^2(0,t;W)}^2+\frac12\Vert w^1-w^2\Vert_{\LDtH}^2+\Vert1* (w^1-w^2)(t)\Vert_{V}^2
\\
\no
&\displaystyle +\sum_{j=1}^3\left(\Vert(\beta_j^1-\beta_j^2)(t)\Vert_{V}^2+\Vert(\beta_j^1-\beta_j^2)\Vert_{\LDtV}^2\right)\\
\no
&
\displaystyle \leq c\itt\sum_{j=1}^3\left(1+\|\dive \ub^1\|_V^2+\|\dive \ub^2\|_V^2\right)\Vert\beta_j^1-\beta_j^2\Vert_{V}^2
\end{eqnarray}
for some positive constant $c$ depending on the data of the problem.
Let us notice that we have estimated the terms containing the nonlinearity $\tau$ in \eqref{p1} and \eqref{p4}
on the right hand side as follows
\begin{eqnarray}\no
&\displaystyle \int_Q\left((\beta_1^1-\beta_2^1)\tau(\teta^1)-(\beta_1^2-\beta_2^2)\tau(\teta^2)\right)\left(\dive (\ub^1-\ub^2)\right)\\
\no
&\displaystyle +\int_Q\left(\left(\tau(\teta^1)-\tau(\teta^2)\right)\dive\ub^1
+\tau(\teta^2)(\dive\ub^1-\dive\ub^2)\right)\left(\left(\beta^1_2-\beta^2_2\right)
-\left(\beta_1^1-\beta_1^2\right)\right)\\
\no
&\displaystyle \leq \frac14 \|w^1-w^2\Vert_{\LDtH}^2+\frac12\Vert (\ub^1-\ub^2)\Vert_{L^2(0,t;W)}^2\\
\no
&\displaystyle + c\itt \left(1+\|\dive\ub^1\|_V^2+\|\dive\ub^2\|_V^2\right)\left(\|\beta_1^1-\beta_2^1\|_V^2+\|\beta_1^2
-\beta_2^2\|_V^2\right)
\,.
\end{eqnarray}
The application of the standard Gronwall lemma together with the regularity \eqref{reg1}
leads to uniqueness of solutions to {\sc Problem (P)} and concludes to proof of
Theorem~\ref{exiteo}.\fin


\section{Proof of Theorem~\ref{uniteo}}
\label{prouni}

In this section we give the proof of Theorem~\ref{uniteo}.  We will use here the same symbol $c$ for some positive constants
(depending only on the data of the problem), which may also be different from line to line.

Then, let us take
two sets of data $(\ub^i_0,\,w_0^i,\,\betab^{i}_0)$,
$({\cal R}^i, \,{\cal F}^i)$ ($i=1,2$) of {\sl Problem (P)} and let $(\ub^i,\,w^i,\,\beta_1^i,\,\beta_2^i,\,\beta_3^i, \,p^i)$
($i=1,2$) be two solutions of {\sl Problem (P)} corresponding to these data.

 Then, write two times equations (\ref{p1}--\ref{p4}) with $(\ub^i,\,w^i,\,\beta_1^i,\,\beta_2^i,\,\beta_3^i,\,p^i)$, make the difference
between the two equations \eqref{p3}, and
test the result with $2(w^1-w^2)$.  Make the difference between the two equations \eqref{p1},
 test the result with $2(\ub^1-\ub^2)_t$. Make the difference between the two equations \eqref{p4},
written for $i=1$ and $i=2$, and test the resulting vectorial equation
by the vector $\big((\beta_1^1)_t-(\beta_1^2)_t, (\beta_2^1)_t-(\beta_2^2)_t, (\beta_3^1)_t-(\beta_3^2)_t\big).$

Finally, summing up the three resulting equations, integrating
 over $(0,t)$, with $t\in [0,T]$, and exploiting
the Lipschitz continuity of $\alpha$ (cf.~assumption \eqref{ipoalfalip}),
we get the following inequality
\begin{eqnarray}\no
&\displaystyle \Vert \ub^1(t)-\ub^2(t)\Vert_W^2+\Vert w^1(t)
-w^2(t)\Vert_H^2+\Vert w^1-w^2\Vert_{\LDtV}^2\\
\no
&\displaystyle +\sum_{j=1}^3\left(\Vert(\beta_j^1)_t-(\beta_j^2)_t\Vert_{\LDtV}^2
+1/2\Vert\nabla(\beta_j^1-\beta_j^2)(t)\Vert_H^2\right)\\
\no
&
\displaystyle \leq c\itt\left(1+\|\dive\ub_1\|_V^2+\|\dive\ub_2\|_V^2\right)\|\teta^1-\teta^2\Vert_{H}^2\\
\no
&\displaystyle +\sum_{j=1}^3\left(c\Vert\beta_j^1-\beta_j^2\Vert_{\LDtH}^2
+1/2\Vert(\beta_j^1)_t-(\beta_j^2)_t\Vert_{\LDtV}^2\right)
\\
\no
&\displaystyle +\|\ub^1_0-\ub_0^2\|_W^2+\| w_0^1-w_0^2\|_H^2+1/2\sum_{j=1}^3\|\nabla(\beta_j^{01}-\beta_j^{02})\|_H^2\\
\label{inequ}
&\displaystyle +c\|{\cal F}^1-{\cal F}^2\|_{L^2(0,T;W')}^2+c\|{\cal R}^1-{\cal R}^2\|_{\LDtVp}^2.
\end{eqnarray}
Let us notice that we have estimated the terms containing the nonlinearity $\tau$ in \eqref{p1} and \eqref{p4}
on the right hand side using \eqref{p2} and the fact that $\gamma$, defined in \eqref{gammateta}, is a locally Lipschitz continuous function,
$\teta^i=\gamma(w^i)$ ($i=1,2$), and $w^i$  are bounded in $L^\infty(Q)$ (cf.~\eqref{reg2}), as follows
\begin{align}\no
&\int_Q\left((\beta_1^1-\beta_2^1)\tau(\teta^1)-(\beta_1^2-\beta_2^2)\tau(\teta^2)\right)\left(\dive (\ub^1-\ub^2)_t\right)\\
\no
&+\int_Q\left(\left(\tau(\teta^1)-\tau(\teta^2)\right)\dive\ub^1+\tau(\teta^2)(\dive\ub^1-\dive\ub^2)\right)\left(\left(\beta^1_2-\beta^2_2\right)_t
-\left(\beta_1^1-\beta_1^2\right)_t\right)\\
\no
&=-\int_Q\left((\beta_1^1-\beta_2^1)\left(\tau(\teta^1)-\tau(\teta^2)\right)+\left((\beta_1^1-\beta_1^2)
-(\beta_2^1-\beta_2^2)\right)\tau(\teta^2)\right)\sum_{j=1}^3\left(\beta_j^1-\beta_j^2\right)_t\\
\no
&\quad+\int_Q\left(\left(\tau(\teta^1)-\tau(\teta^2)\right)\dive\ub^1
+\tau(\teta^2)(\dive\ub^1-\dive\ub^2)\right)\left(\left(\beta^1_2-\beta^2_2\right)_t
-\left(\beta_1^1-\beta_1^2\right)_t\right)\\
\no
&\leq \frac14\sum_{j=1}^3\|(\beta_j^1-\beta_j^2)_t\|_V^2+c\sum_{j=1}^3\|\beta_j^1-\beta_j^2\|_{\LDtH}^2\\
\no
&\quad+ c\itt \left(1+\|\dive\ub^1\|_V^2+\|\dive\ub^2\|_V^2\right)\left(\|w^1-w^2\Vert_{H}^2+\Vert (\ub^1-\ub^2)\Vert_{W}^2\right)
\,.
\end{align}
Moreover, by adding to both sides in the inequality \eqref{inequ}
\begin{equation}\no
1/2\sum_{j=1}^3\Vert(\beta_j^1)(t)-(\beta_j^2)(t)\Vert_{H}^2=1/2\sum_{j=1}^3\|\beta_j^{01}-\beta_j^{02}\|_H^2
+\itt\sum_{j=1}^3((\beta_j^1)_t-(\beta_j^2)_t,\beta_j^1-\beta_j^2)
\end{equation}
 we get the following inequality
\begin{eqnarray}\no
&\displaystyle \Vert \ub^1(t)-\ub^2(t)\Vert_W^2+\Vert w^1(t)-w^2(t)\Vert_H^2+\Vert w^1-w^2\Vert_{\LDtV}^2\\
\no
&\displaystyle +\sum_{j=1}^3\left(\Vert(\beta_j^1)_t-(\beta_j^2)_t\Vert_{\LDtV}^2+\Vert(\beta_j^1-\beta_j^2)(t)\Vert_V^2\right)\\
\no
&\displaystyle
\leq c\itt\left(1+\|\dive\ub_1\|_V^2+\|\dive\ub_2\|_V^2\right)\left(\|w^1-w^2\Vert_{H}^2+\Vert (\ub^1-\ub^2)\Vert_{W}^2\right)\\
\no
&\displaystyle +c\Big(\sum_{j=1}^3\Vert\beta_j^1-\beta_j^2\Vert_{\LDtH}^2+\|\ub^1_0-\ub_0^2\|_W^2+\| w_0^1-w_0^2\|_H^2
+\sum_{j=1}^3\|\beta_j^{01}-\beta_j^{02}\|_V^2
\\
\no
&\displaystyle+\|{\cal F}^1-{\cal F}^2\|_{L^2(0,T;W')}^2+\|{\cal R}^1-{\cal R}^2\|_{\LDtVp}^2\Big)\,.
\end{eqnarray}
Applying now a standard version of Gronwall's  lemma (cf.~\cite[Lemme~A.4, p.~156]{brezis}),
we get the desired continuous dependence estimate \eqref{CD}.
This concludes the proof of Theorem~\ref{uniteo}.\fin


\end{document}